\documentclass[reqno]{amsart}
\usepackage[utf8]{inputenc}
\usepackage{fullpage,hyperref,enumitem}
\usepackage{tikz,amsmath,amssymb,color}
\usepackage{graphicx}
\usepackage[margin=1in]{geometry}
\usetikzlibrary{calc}
\usepackage{float}
\usepackage{amsrefs}
\usepackage[mathscr]{euscript}
\usepackage{colonequals}
\usepackage{stackengine}
\usepackage{xcolor}
\usepackage{mathtools}

\newcommand{\rthmname}{Theorem}
\newtheorem*{rthm}{\rthmname}
\newenvironment{restatethm}[1]{\renewcommand{\rthmname}{#1}
\begin{rthm}
}{\end{rthm}}
\theoremstyle{definition} 
\newtheorem{theorem}{Theorem}[section]
\newtheorem{corollary}{Corollary}[theorem]
\newtheorem{proposition}[theorem]{Proposition}

\newtheorem{example}[theorem]{Example}

\newtheorem{lemma}[theorem]{Lemma}

\newtheorem{definition}[theorem]{Definition}
\newtheorem*{block}{Block Decomposition}
\theoremstyle{remark}

\newcommand{\NN}{\mathbb{N}}

\newcommand{\id}{\mathrm{id}}

\newcommand{\PF}{\mathrm{PF}}

\newcommand{\out}{\mathcal{O}}
\newcommand{\Lucky}{\mathsf{Lucky}}
\newcommand{\lucky}{\mathsf{lucky}}

\newcommand{\Sym}{\mathfrak{S}}

\DeclareMathOperator\Av{Av}

\newcommand{\delete}[1]{}

\newcommand{\seqnum}[1]{\href{https://oeis.org/#1}{\rm  \underline{#1}}}
\usepackage{thmtools} 
\usepackage{thm-restate}

\title{Pattern avoidance classes of parking permutations with fixed lucky sets}

\author[Ferreri]{Melanie Ferreri}
\address[M.~Ferreri]{Department of Mathematics, William \& Mary, Williamsburg, VA 23187}
\email{\textcolor{blue}{\href{mailto:email}{mjferreri@wm.edu}}}

\author[Martinez]{Lucy Martinez}
\address[L.~Martinez]{Department of Mathematics, Rutgers University, Piscataway, NJ 08854}
\email{\textcolor{blue}{\href{mailto:lucy.martinez@rutgers.edu}{lucy.martinez@rutgers.edu}}}

\begin{document}

\begin{abstract}
For a parking function, the resulting sequence after all cars have parked is called the associated parking permutation. For each parking function, cars that park in their preferred spot are called lucky cars.
In this paper, we enumerate parking permutations that avoid a specified set of permutation patterns in $\Sym_3$ and whose lucky cars are the first $k$ cars. Often, the enumeration simplifies to a product of a fraction with a binomial coefficient, recovering existing sequences in the OEIS. For two of the classes, the enumerations are new sequences.
\end{abstract}

\maketitle

\section{Introduction}
\emph{Parking functions} have been central to many active areas of research since they were implicitly studied by Pyke~\cite{Pyke} and later on by Konheim and Weiss~\cite{konheim1966occupancy}. As integer sequences associated to permutations, questions of statistics and pattern avoidance for these sequences have also been popular in the literature. One such statistic is the $\lucky$ statistic for parking functions.
This statistic is closely related to descents and ties, and these are of interest partially due to their algebraic and geometric connections; for example, descents in tieless parking functions have been studied in connection with the Chow polynomial of the noncrossing partition lattice~\cites{chun2026chowpolynomialnoncrossingpartition,alexandersson2026parkingfunctionssmirnovwords}, and for vector parking functions with a specified parking protocol, Avila, Ferroni, and Morales recently demonstrated that the $\lucky$ statistic appears in a combinatorial interpretation of coefficients for the Ehrhart polynomial of the Pitman-Stanley polytope corresponding to a vector~\cite[Theorem 1.1]{luckandmagic}. For other statistics on parking functions, see for example~\cites{statsinPFs, Schumacher}, and for a survey on parking functions, see Yan~\cite{Yan2015}.

The study of the $\lucky$ statistic for parking functions began with the work of Gessel and Seo, who provided a polynomial expression for the generating function that enumerates parking functions with a fixed number of lucky cars~\cite[Theorem 10.1]{GesselSeo}. Later, Shin provided a bijection that maps labeled forests with a certain number of leaders to parking functions with the same number of lucky cars~\cite{shin2008newbijectionforestsparking}.
Harris and Martinez gave a characterization for the set of parking permutations arising from parking functions with a fixed set of lucky cars~\cite{harrismartinez}. More recently, Ferreri, Harris, Martinez, and Swartz proved that every parking permutation admits a decomposition into $|I|$ increasing subsequences, where each subsequence begins with an element of the lucky set $I$~\cite{vectorparkingfunctionsluckysets}.
The $\lucky$ statistic was also considered by Colmenarejo et al.~\cite{colmenarejoluckystirling}, who studied Stirling permutations as parking functions and considered the question of how many Stirling permutations have certain lucky sets. For more work on the $\lucky$ statistic and its connections, see~\cites{luckycarsandspots, quicksort, luckycarsexpectedvalue}.

Pattern avoidance in parking functions has also been studied in different contexts. Motivated by Remmel and Qiu~\cite{QiuRemmel}, Adeniran and Pudwell studied pattern avoidance in parking functions by viewing parking functions as labeled Dyck paths~\cite{AdeniranPudwell}. They enumerated parking functions avoiding any set of two or more patterns of length $3$, and Adenbaum later completed the enumeration for all single patterns of length $3$~\cite{Adenbaum26}.
Another notion of pattern avoidance in parking functions was studied by Yan~\cite{JunYan}. In this version, Yan studied the parking permutation map, which is the order in which the cars ultimately park on the street, and enumerated all parking functions whose parking permutations avoid any collection of patterns of length three.

Motivated by recent results on pattern avoidance and the $\lucky$ statistic for parking functions, we investigate the following question: given a fixed lucky set $I=[k]$, what structure is imposed on parking permutations with lucky set $I$ that avoid a prescribed pattern? In the present work, we provide a solution to this question for classical patterns of length $3$, and show that the resulting parking permutations admit a structure in terms of ascent and descent statistics of the permutation induced by the lucky cars. Our main enumeration results are as follows. For four of the patterns, the enumeration matched existing sequences in the OEIS \cite{OEIS}.
\begin{restatethm}{Theorem~\ref{thm:123avoid}}
[Lobb numbers, OEIS \seqnum{A039599}]
The number of parking permutations of length $n$ that avoid the pattern $123$ with lucky set $I=\{1,2,\ldots,k\}$ is given by
\[\frac{2n-2k+1}{n+1}\binom{2k}{n}.
\]
\end{restatethm}

\begin{restatethm}{Theorem~\ref{thm:321avoid}}
Let $k\geq 2$. The number of parking permutations of length $n$ that avoid the pattern $321$ with lucky set $I=\{1,2,\ldots,k\}$ is given by
\[\sum_{i=0}^{n-k} \binom{n-k}{i}\left(\frac{n-i-k+2}{n-i+1}\binom{n-i+k-1}{k-1} -\binom{n-i-2}{k-1}\right). 
\]
\end{restatethm}

\begin{restatethm}{Theorems~\ref{thm:213avoid}, \ref{thm:231avoid}, \ref{thm:132avoid}}
[Catalan's triangle (read by rows), OEIS \seqnum{A009766}]
The number of parking permutations of length $n$ that avoid the pattern $213$, or the pattern $231$, or the pattern $132$, with lucky set $I=\{1,2,\ldots, k\}$ is given by
\[\frac{n-k+2}{n+1}\binom{n+k-1}{k-1}.
\]
\end{restatethm}

\begin{restatethm}{Theorem~\ref{thm:312avoid}}
The number of parking permutations of length $n$ that avoid the pattern $312$ with lucky set $I=\{1,2,\ldots, k\}$ is given by
\[\frac{n-k+1}{n+1}\binom{n+k}{k} -\binom{n-1}{k}.\]
\end{restatethm}

Tables of values for each of these formulas are given in Appendix~\ref{appendix}.

\section{Definitions and notation}
Let $n\in \NN=\{1,2,\ldots\}$. A \emph{parking function} of length $n$ is an integer sequence $\alpha=(a_1, a_2,\ldots, a_n)$ whose nondecreasing rearrangement $(a_{(1)}, a_{(2)},\ldots, a_{(n)})$ satisfies $a_{(i)} \leq i$ for each $1\leq i \leq n$. 
Such sequences are named for their interpretation in terms of cars parking in a one-way street: the entry $a_i$ indicates the preferred parking spot of the $i$th car. Each car first attempts to park in its preferred spot, and if that spot is not available, checks each subsequent spot in order until it finds an available spot, and parks there. The sequence is considered a parking function only if all cars are able to park under this protocol. We let $\PF_n$ denote the set of parking functions of length $n$.

Once all $n$ cars have parked, the resulting arrangement of cars is called the \emph{outcome}, or \emph{parking permutation} corresponding to the parking function $\alpha$. Let $[n] = \{1,2,\ldots, n\}$ and let $\Sym_n$ denote the set of permutations of $[n]$. Then the parking permutation of a parking function $\alpha$ is denoted by \[\out(\alpha)=\pi_1\pi_2\cdots \pi_n,\]
where $\pi_i$ is the car that parked in spot $i$.
If a car parks in its preferred spot, that car is said to be \emph{lucky}; otherwise it is \emph{unlucky}. For a given parking function $\alpha$, we denote by $\Lucky(\alpha)$ the set of all lucky cars in $\alpha$, and define $\lucky(\alpha) = |\Lucky(\alpha)|$. 

A permutation $\sigma \in \Sym_n$ is said to have an \emph{ascent} at position $i$ if $\sigma_i<\sigma_{i+1}$, and $\sigma$ has a \emph{descent} at position $i$ if $\sigma_i>\sigma_{i+1}$.
Harris and Martinez proved that in the parking permutation of a given parking function, the location of the descents force certain cars to be in the lucky set~\cite[Lemma 2.2]{harrismartinez}. Moreover, the first car to prefer spot $1$ will park there; hence that car is a lucky car~\cite[Theorem 2.3]{harrismartinez}. These facts are key for our results, and for convenience, we rephrase them below.

\begin{lemma}[\cite{harrismartinez}, Lemma 2.2]
\label{lem:luckydescents}
Fix a lucky set $I\subseteq [n]$ of $\PF_n$. Let $\pi\in \Sym_n$ be the parking permutation of a parking function with lucky set $I$. If $\pi$ has a descent at position $j$, then $\pi_{j+1} \in I$. Consequently, if $u\notin I$ then $u$ cannot be directly preceded by any $v>u$.
\end{lemma}

\begin{lemma}[\cite{harrismartinez}, Theorem 2.3]\label{lem:firstcarislucky}
Fix a lucky set $I\subseteq [n]$ of $\PF_n$. If $\pi\in \Sym_n$ is the parking permutation of a parking function of length $n$ with lucky set $I$, then $\pi_1\in I$.
\end{lemma}

We denote by $\out_n(I)$ the set of parking permutations of length $n$ corresponding to any parking function with lucky set $I$. That is,
\[
\out_n(I) = \{\sigma \in \Sym_n\, : \sigma = \out(\alpha) \text{ where }\Lucky(\alpha) = I\}.
\]
Throughout this paper, we restrict our results to the case where $I=\{1,2,\ldots,k\}$, with $k\leq n$. 

Let $\sigma \in \Sym_m$. A permutation $\pi\in \Sym_n$ is said to \emph{contain} $\sigma$ as a pattern if there is a length-$m$ subsequence of $\pi$ whose entries appear in the same relative order as the entries in $\sigma$. Otherwise, $\pi$ is said to \emph{avoid} $\sigma$. Let $\Pi=\{\sigma_1,\sigma_2,\ldots, \sigma_i\}$ be a collection of patterns and let $\Av_n(\Pi)$ be the set of permutations that avoid all patterns in the set $\Pi$. 
For more background on pattern avoidance, see B\'ona~\cite[Chapter 4]{MBona12}.

Let $\out_n(I; \, \sigma_1,\sigma_2,\ldots,\sigma_i)$ denote the set of parking permutations with lucky set $I$ that avoid each of the patterns $\sigma_1,\sigma_2,\ldots,\sigma_i$. 
Although there are Catalan many permutations that avoid a single pattern of length $3$, not every such permutation arises as a parking permutation with a fixed set of lucky cars. This is illustrated in the following example.

\begin{example}
    Let $n=5$ and $I=\{1,2,3\}$. We check that $\out_5(I;\, 123) \subsetneq \Av_5(123)$.
    Indeed, consider $\pi=35412\in \Av_5(123)$. Since there is a descent at position $i=2$, by Lemma~\ref{lem:luckydescents} car $4$ must belong to the set $I$. Hence, $35412\notin \out_5(I; \, 123)$.
\end{example}

In some of our results, we make use of the following result by Connolly, Gabor, and Godbole, who considered the question of how many permutations avoid the pattern $123$ and whose first ascent occurs at position $k$.
\begin{theorem}[Connolly, Gabor, and Godbole~\cite{firstascent}]\label{thm:ConnollyGaborGodbole}
The number of $123$-avoiding permutations of length $n$ whose first ascent occurs at position $k$ is given by the $k$-fold Catalan convolution $C_{n,k}$, where
\begin{equation*}
C_{n,k}=\sum_{\substack{i_1+\cdots+i_k=n, \\i_j\geq 1}}\,\prod_{r=1}^k C_{i_r-1}=\frac{k}{2n-k}\binom{2n-k}{n}, \qquad C_n=\frac{1}{n+1}\binom{2n}{n}.
\end{equation*}
\end{theorem}

\section{Parking permutation decomposition avoiding a single pattern of length three}

The following alternative definition of $\out_n(I; \, \sigma_1,\sigma_2,\ldots,\sigma_i)$ gives a characterization of the elements in this set. The first two criteria are shown to be necessary in~\cite[Theorem 2.3]{harrismartinez}.

\begin{definition}
\label{def:characterization}
    Fix a lucky set $I \subseteq [n]$. Then 
    $\pi = \pi_1 \pi_2 \cdots \pi_n \in \Sym_n$ is in $\out_n(I;\, \sigma_1, \sigma_2,\ldots, \sigma_i)$ if and only if
    \begin{enumerate}
    \item $\pi_1 \in I$,
    \item if $\pi_{i-1} > \pi_i$, then $\pi_i \in I$, and
    \item $\pi$ avoids each of the patterns $\sigma_1, \sigma_2,\ldots, \sigma_i$.
    \end{enumerate}
\end{definition}

\begin{example}\label{ex:avoiding perms ex}
Let $n=5$ and $I=\{1,2,3\}$.
If $\sigma=123$ then 
\[\out_5(I; \, 123)=\{ 25143, 32514, 35142, 35214, 35241\}.\]
Note that there are $42$ permutations that avoid the pattern $123$ of length $5$, while the set of $123$-avoiding parking permutations with lucky set $I$ contains only $5$ elements.

As a different example, if $\sigma=321$, then the set $\out_5(I;\, 321)$ contains the following $32$ elements,
\[ 
\out_5(\{1,2,3\}; \, 321)=
    \left\{
    \begin{array}{cccccccc}
    12345 & 12435 & 12453 & 12534 & 13245 & 13425 & 13452 & 13524 \\
    14235 & 14253 & 14523 & 15234 & 21345 & 21435 & 21453 & 21534 \\
    23145 & 23415 & 23451 & 23514 & 24135 & 24153 & 24513 & 25134 \\
    31245 & 31425 & 31452 & 31524 & 34125 & 34152 & 34512 & 35124
    \end{array}
    \right\}. \]
\end{example}

\subsection{The class $\Av(123)$ with fixed lucky sets}\label{subsec:123avoid}

We first provide a characterization of $123$-avoiding parking permutations with lucky set $[k]$. This construction is implicit in the proof of Theorem~\ref{thm:123avoid}.

\begin{block}
Let $I=\{1,2,\ldots,k\}\subseteq [n]$. For each permutation $\pi=\pi_1\pi_2\cdots\pi_n\in \out_n(I;\, 123)$, there exists a permutation
\[\sigma=\sigma_1\sigma_2\cdots\sigma_k \in \Av_k(123),\]
whose first ascent occurs at position $i$ ($\sigma_i<\sigma_{i+1})$, taking $i=k$ if $\sigma$ has no ascents, such that $\pi$ admits a block decomposition
\[\pi=B_1B_2\cdots B_k\in \out_n(I;\, 123),\]
satisfying:
\begin{enumerate}
\item For any $r\in [k]$, the first element of the block $B_r$ is the lucky car $\sigma_r$.
\item If $r>i$, then $B_r=\sigma_r$.
\item The subsequence $\sigma_1\sigma_2\cdots\sigma_i$ is a decreasing subsequence since $i$ is the position where the first ascent occurs.
\item For $r\leq i$, each block $B_r$ consists of the lucky car $\sigma_r$ followed by at most one unlucky car from $[n]\setminus I=\{k+1,k+2,\ldots,n\}$.
\item For any $u,v\in [n]\setminus [k]$, if $u$ and $v$ are unlucky cars belonging to blocks $B_r$ and $B_s$, respectively, with $r<s\leq k$, then $u>v$.
\end{enumerate}
\end{block}

\begin{example}
Let $n=9$ and $k=6$. The following are examples of parking permutations in $\Av_9(123)$ with lucky set $[6]$. The blocks are indicated with underbraces.

\[
   \pi= \underbrace{69}_{B_1} \underbrace{58}_{B_2} \underbrace{4}_{B_3} \underbrace{2}_{B_4} \underbrace{17}_{B_5}
   \underbrace{3}_{B_6}
\]

\[ \sigma = \underbrace{5}_{B_1} \underbrace{49}_{B_2} \underbrace{38}_{B_3} \underbrace{27}_{B_4} \underbrace{6}_{B_5} \underbrace{1}_{B_6} \]

\[ \tau = \underbrace{59}_{B_1} \underbrace{48}_{B_2} \underbrace{37}_{B_3} \underbrace{6}_{B_4} \underbrace{1}_{B_5}\underbrace{2}_{B_6}\]
\end{example}

\begin{lemma}
\label{lem:123-avoiding decreasing unlucky}
In a $123$-avoiding parking permutation with lucky set $[k]$, the unlucky cars $\{k+1,k+2,\ldots, n\}$ must appear in decreasing order.
\end{lemma}

\begin{proof}
Assume for contradiction that we have two unlucky cars $u_1, u_2\in \{k+1,k+2,\ldots, n\}$ such that $u_1<u_2$ and $u_2$ appears to the right of $u_1$. Since $u_1$ is unlucky, $u_1$ must appear after a lucky car $\sigma_a$, and so in the parking permutation we have the subsequence $\sigma_a,u_1,u_2$. Observe that $u_1,u_2\geq k+1$ and $\sigma_a \in [k]$. Therefore, $\sigma_a<u_1<u_2$ forms an increasing subsequence in the parking permutation, which is a contradiction. Hence, for any unlucky cars $u_1,u_2\in U, u_1<u_2$, the unlucky car $u_2$ must appear before the unlucky car $u_1$.
\end{proof}

\begin{lemma}\label{lem:123avoiding}
In a $123$-avoiding parking permutation with lucky cars $\{1,2,\ldots, k\}$, let $\sigma = \sigma_1 \sigma_2 \cdots \sigma_k\in \Sym_k$ be the subsequence consisting of all lucky cars. If there is an ascent at position $i$ in $\sigma$, then all unlucky cars must appear in positions to the left of $\sigma_{i+1}$.
\end{lemma}

\begin{proof}
Suppose $\sigma$ has an ascent at position $i$, so $\sigma_{i}<\sigma_{i+1}$. Suppose to the contrary that there is some unlucky car $\tau$ appearing to the right of $\sigma_{i+1}$. Then, since $\sigma_{i}<\sigma_{i+1}\leq k$ and $\tau\geq k$ as $\tau\in [n]\setminus[k]$, we have that $\sigma_{i},\sigma_{i+1},\tau$ forms an increasing subsequence in the parking permutation, which is a contradiction since the permutation must avoid $123$.
\end{proof}

\begin{theorem}[Lobb numbers, OEIS \seqnum{A039599}]\label{thm:123avoid}
The number of parking permutations of length $n$ that avoid the pattern $123$ with lucky set $I=\{1,2,\ldots,k\}$ is
\[ |\out_n(I;\, 123)| =\frac{2n-2k+1}{n+1}\binom{2k}{n}.
\]
\end{theorem}
\begin{proof}
We count the number of ways to construct a parking permutation that avoids $123$.
First, we choose a $123$-avoiding permutation $\sigma$ of $[k]$; this determines the order in which the lucky cars $\{1,2,\ldots, k\}$ appear in the parking permutation. 

Then we place the remaining $n-k$ unlucky cars. Each unlucky car must be placed after a lucky car, and by Lemma \ref{lem:luckydescents}, for sequences of consecutive unlucky cars appearing between lucky cars, the unlucky cars must be increasing. However, since the parking permutation must avoid $123$, these consecutive subsequences can never have length greater than 2; that is, each lucky car can be followed by at most one unlucky car. So for each unlucky car, we choose a (unique) lucky car for it to be placed directly after.

Moreover, when placing the unlucky cars, we must ensure that the resulting parking permutation still avoids $123$. 
If $\sigma$ has an ascent at $\ell$, then by Lemma~\ref{lem:123avoiding}, we cannot place any unlucky cars after $\sigma_{\ell+1}$. Let $i$ be the position of the first ascent of $\sigma$ and if $\sigma$ has no ascents, we let $i=k$. Then all unlucky cars must be placed before position $i+1$.

So we may insert each of the $n-k$ unlucky cars after any of the lucky cars at positions $\{1,2,\ldots, i\}$ (again, with at most one unlucky car after each lucky car). There are then $\binom{i}{n-k}$ choices for how to place the unlucky cars. Also, once we choose where to place them, by Lemma \ref{lem:123-avoiding decreasing unlucky}, we must insert them in decreasing order (so their order is fixed).

So for each $123$-avoiding permutation of $[k]$ which has first ascent at position $i$, there are $\binom{i}{n-k}$ ways to construct a parking permutation of length $n$.
By Theorem~\ref{thm:ConnollyGaborGodbole}, the number of $123$-avoiding permutations of $[k]$ which have first ascent at position $i$ is given by the $i$-fold Catalan convolution formula $C_{k,i} = \frac{i}{2k-i}\binom{2k-i}{k}$. So in total, the number of possible parking permutations is $\sum_{i=1}^k \binom{i}{n-k}C_{k,i}$.

The position $i$ of the first ascent can range from $1$ to $k$, but if $i<n-k$, then $\binom{i}{n-k} = 0$. So this formula may be rewritten as
\[
\sum_{i=1}^k \binom{i}{n-k}C_{k,i} = \sum_{i=n-k}^k \binom{i}{n-k}C_{k,i}=\frac{2n-2k+1}{n+1}\binom{2k}{n},
\]
where the last equality is proven in Lemma~\ref{lem:identity1}
.\qedhere
\end{proof}

\subsection{The class $\Av(321)$ with fixed lucky sets}
The characterization of the $123$-avoiding parking permutations provides insight to the structure of the parking permutations that avoid $321$ via the reversal map.

Every parking permutation in $\out_n([k];\, 321)$ admits a block decomposition as follows. One can check that such conditions are sufficient for a permutation to satisfy Definition~\ref{def:characterization}, and they are shown to be necessary in the remaining proofs in this section.

\begin{block}
Let $I=\{1,2,\ldots,k\}\subseteq [n]$. For each permutation $\pi=\pi_1\pi_2\cdots \pi_n \in \out_n(I;\, 321)$, there exists a permutation
\[\sigma=\sigma_1\sigma_2\cdots\sigma_k \in \Av_k(321),\]
whose last descent occurs at position $j$ ($\sigma_j>\sigma_{j+1}$), taking $j=0$ if $\sigma$ has no descents, such that $\pi$ admits a block decomposition
\[\pi=B_1B_2\cdots B_{k} \in \out_n(I; \, 321),\]
satisfying:
\begin{enumerate}
    \item For any $r\in [k]$, the first element of the block $B_r$ is the lucky car $\sigma_r$.
    \item If $r<j$ then $B_r=\sigma_r$.
    \item The subsequence $\sigma_{j+1}\sigma_{j+2}\cdots\sigma_k$ is an increasing subsequence since $j$ is the position where the last descent occurs.
    \item For $r\geq j$, each block $B_r$ consists of the lucky car $\sigma_r$ followed by an increasing sequence (possibly empty) of unlucky cars from $[n]\setminus I=\{k+1,k+2,\ldots,n\}$. 
    \item For any $u,v\in [n]\setminus [k]$, if $u$ and $v$ are unlucky cars belonging to blocks $B_r$ and $B_s$, respectively, with $r<s<k$, then $u<v$.
\end{enumerate}
\end{block}

\begin{example}
Let $n=9$ and $k=5$. The following are examples of parking permutations in $\Av_9(321)$ with lucky set $[5]$. The blocks are indicated with underbraces and each element at the beginning of each block is a lucky car.

\[
   \pi= \underbrace{4}_{B_1} \underbrace{1}_{B_2} \underbrace{5}_{B_3} \underbrace{267}_{B_4} \underbrace{389}_{B_5}
\]

\[ \sigma = \underbrace{3}_{B_1} \underbrace{17}_{B_2} \underbrace{2}_{B_3} \underbrace{49}_{B_4} \underbrace{568}_{B_5}\]

\[ \tau = \underbrace{2}_{B_1} \underbrace{1}_{B_2} \underbrace{3}_{B_3} \underbrace{56789}_{B_4} \underbrace{4}_{B_5}\]
\end{example}

\begin{definition}[Reverse Permutation]
Let $\sigma=\sigma_1\sigma_2\cdots\sigma_n\in \Sym_n$. We define the reverse permutation $\sigma^r=\sigma_n\sigma_{n-1}\cdots\sigma_2\sigma_1\in \Sym_n$.
\end{definition}

There is a natural bijection between the $321$-avoiding permutations of length $k$ whose last descent appears at some position $j$ and the $123$-avoiding permutations of the same length whose first ascent is at position $k-j$ via the reversal map. 

\begin{lemma}[Correspondence between $123$- and $321$-avoiding permutations]\label{lemma:reverse bijection}
Let $\pi=\pi_1\pi_2\cdots\pi_k\in \Av_k(321)$, with last descent at position $j$ ($\pi_j>\pi_{j+1}$). Then, $\pi^r=\pi_k\pi_{k-1}\cdots\pi_1\in \Av_{k}(123)$ and $\pi^r$ has first ascent at position $k-j$ ($\pi_{k-j}<\pi_{k-j+1}$).
\end{lemma}

\begin{proof}
It is a classical result that the pattern $123$ is Wilf-equivalent to the pattern $321$, because $321$ is the reverse of $123$, and the reverse of every permutation that avoids $123$ avoids $321$. Consider $\pi\in \Av_k(321)$ such that the last descent of the permutation $\pi$ appears at position $j$. Then, by reversing this permutation $\pi$, we obtain 
\[\pi^r=\pi_k\pi_{k-1}\cdots \pi_1,\]
such that $\pi^r_i=\pi_{k-i+1}$. It follows that
\[\pi^r_{k-j} = \pi_{k-(k-j)+1} = \pi_{j+1}, \]
and
\[\pi^r_{k-j+1} = \pi_{k-(k-j+1)+1} = \pi_{j}>\pi_{j+1} = \pi^r_{k-j}.\]
Therefore, $\pi^r_{k-j+1}>\pi^r_{k-j}$, and so $\pi^r$ has an ascent at position $k-j$. Moreover, for all $m<k-j$,
\[\pi^r_m=\pi_{k-m+1}<\pi_{k-m+2}=\pi^r_{m-1}.\]
Hence, $\pi^r$ has no ascents before position $k-j$. Thus, $\pi^r$ is a $123$-avoiding permutation whose first ascent is at position $k-j$.
\end{proof}

\begin{corollary}(Catalan's triangle (read by rows), OEIS \seqnum{A009766})
\label{cor:enumerationof321avoidingpermutationsbylastdescent}
    The number of $\sigma \in \Av_k(321)$ whose last descent appears at position $j$ is given by $C_{k,k-j} = \frac{k-j}{k+j}\binom{k+j}{k}$. 
\end{corollary}

\begin{proof}
    By Lemma \ref{lemma:reverse bijection}, such permutations are in bijection with the set of $\sigma \in \Av_k(123)$ whose first ascent appears at position $k-j$. By Theorem~\ref{thm:ConnollyGaborGodbole}, it follows that these permutations are enumerated by the $(k-j)$-fold Catalan convolution formula \[C_{k,k-j} = \frac{k-j}{2k-(k-j)}\binom{2k-(k-j)}{k} = \frac{k-j}{k+j}\binom{k+j}{k}. \qedhere\]
\end{proof}

\begin{lemma}
\label{lem:321avoiding1}
In a $321$-avoiding parking permutation with lucky cars $\{1,2,\ldots, k\}$, let $\sigma = \sigma_1 \sigma_2 \cdots \sigma_k\in \Sym_k$ be the subsequence consisting of all lucky cars. For all unlucky cars which appear before $\sigma_k$, they must appear in increasing order.
\end{lemma}

\begin{proof}
Suppose to the contrary that for unlucky cars $\tau_i, \tau_j \in \{k+1, k+2,\ldots, n\}$, we have $\tau_i>\tau_j$ and $\tau_i$ appears before $\tau_j$ in the parking permutation, with both $\tau_i$ and $\tau_j$ appearing before $\sigma_k$. Because both $\tau_i$ and $\tau_j$ are greater than $\sigma_k \in [k]$, it follows that $\tau_i > \tau_j >\sigma_k$ forms a decreasing subsequence in the parking permutation, which is a contradiction. So the subsequence consisting of all unlucky cars appearing before $\sigma_k$ must be strictly increasing.
\end{proof}

\begin{lemma}
\label{lem:321avoiding2}
In a $321$-avoiding parking permutation with lucky cars $\{1,2,\ldots, k\}$, let $\sigma = \sigma_1 \sigma_2 \cdots \sigma_k\in \Sym_k$ be the subsequence consisting of all lucky cars. If there is a descent at position $j$ in $\sigma$, then all unlucky cars must appear in positions to the right of $\sigma_j$.
\end{lemma}

\begin{proof}
Suppose $\sigma$ has a descent at position $j$, so $\sigma_j > \sigma_{j+1}$. Suppose to the contrary that there is some unlucky car $\tau$ appearing to the left of $\sigma_j$. Then since $\sigma_j \leq k$ and $\tau>k$, we have that $\tau>\sigma_j>\sigma_{j+1}$ forms a decreasing subsequence in the parking permutation, which is a contradiction as the permutation must avoid $321$.
\end{proof}

\begin{theorem}
\label{thm:321avoid}
Let $k\geq 2$. The number of parking permutations of length $n$ that avoid the pattern $321$ with lucky set $I=\{1,2,\ldots,k\}$ is
\[|\out_n(I; \, 321)|= \sum_{i=0}^{n-k} \binom{n-k}{i}\left(\frac{n-i-k+2}{n-i+1}\binom{n-i+k-1}{k-1} -\binom{n-i-2}{k-1}\right). 
\]
If $k=1$ and $n\in \NN$ then $|\out_n(I; \, 321)|=1$.
\end{theorem}

\begin{proof}
Suppose $k\geq 2$ and let $\{1,2,\ldots,k\}$ be the set of lucky cars.
To construct a $321$-avoiding parking permutation, we begin by choosing a permutation $\sigma \in \Av_k(321)$ which determines the placement of the lucky cars $\{1,2,\ldots, k\}$. Then, we place the unlucky cars $\{k+1, k+2, \ldots, n\}$.

First fix $i \in \{0,1,\ldots, n-k\}$, and choose $i$ of the unlucky cars to appear after the last lucky car $\sigma_k$. Since these together form a contiguous block of unlucky cars, once they are selected, they are fixed to appear in increasing order by Lemma~\ref{lem:luckydescents}. This choice can be made in $\binom{n-k}{i}$ ways.

Next, the remaining $n-k-i$ unlucky cars are then placed to the left of $\sigma_k$. (Since $k\geq 2$, there is at least one lucky car before $\sigma_k$, so unlucky cars can be placed to its left.)

If $\sigma=12\cdots k=\mathrm{id}$, the remaining unlucky cars can be distributed into any of the ``gaps'' directly following cars $\sigma_1,\sigma_2,\ldots, \sigma_{k-1}$. Distributing the $n-k-i$ unlucky cars over these $k-1$ gaps amounts to choosing a weak composition of $n-k-i$ into $k-1$ parts. By Lemma \ref{lem:321avoiding1}, once their positions are chosen, they are fixed to appear in increasing order. So this can be done in $\binom{n-i-1-1}{k-1-1} = \binom{n-i-2}{k-2}$ ways.

Otherwise, suppose that $\sigma \neq \mathrm{id}$, so there exists some descent in $\sigma$. Let $j$ be the position of the last descent in $\sigma$. 
From Corollary \ref{cor:enumerationof321avoidingpermutationsbylastdescent}, we have that the number of possible $\sigma \in \Av_k(321)$ whose last descent is in position $j$ is given by $C_{k,k-j} = \frac{k-j}{k+j}\binom{k+j}{k}$, and $j$ can take any value from $1$ to $k-1$.
By Lemma \ref{lem:321avoiding2}, the unlucky cars must all appear after $\sigma_j$, and again by Lemma \ref{lem:321avoiding1}, once their positions are chosen, they are fixed to appear in increasing order. So we may distribute these $n-k-i$ cars into the $k-j$ possible positions directly following cars $\sigma_j, \sigma_{j+1},\ldots,\sigma_{k-1}$. This can be done in $\binom{n-k-i+k-j-1}{k-j-1} = \binom{n-i-j-1}{k-j-1}$ ways.

Then in total, the number of $321$-avoiding parking permutations with lucky set $\{1,2,\ldots,k\}$ is 
\begin{align*}
  |\out_n(I; \, 321)| &=  \underbrace{\sum_{i=0}^{n-k}\binom{n-k}{i} \binom{n-i-2}{k-2}}_{\text{if $\sigma=\mathrm{id}$}} + \underbrace{\sum_{i=0}^{n-k}\binom{n-k}{i} \sum_{j=1}^{k-1}  \binom{n-i-j-1}{k-j-1}C_{k,k-j} }_{\text{if $\sigma\neq\mathrm{id}$}} \\
  &=  \sum_{i=0}^{n-k} \binom{n-k}{i}\left(  \binom{n-i -2}{k-2} + \sum_{j=1}^{k-1} \binom{n-i-j-1}{k-j-1}C_{k,k-j}\right) \\
  &= \sum_{i=0}^{n-k} \binom{n-k}{i}\left(  \binom{n-i -2}{k-2} +  \frac{n-i-k+2}{n-i+1}\binom{n-i+k-1}{k-1} - \binom{n-i-1}{k-1} \right) \\
&=\sum_{i=0}^{n-k} \binom{n-k}{i} \left(\frac{n-i-k+2}{n-i+1}\binom{n-i+k-1}{k-1} -\binom{n-i-2}{k-1}\right),
\end{align*}
where the last two equalities follow from Lemma~\ref{lem:321simplify} and Pascal's rule, respectively.

Finally, if $k=1$, then the lucky car $1$ must appear in position $1$, and all remaining unlucky cars must be placed after this in increasing order, yielding the identity as the only possible parking permutation.
\end{proof}

\subsection{The class $\Av(213)$ with fixed lucky sets}

We begin with the block decomposition for the $213$-avoiding parking permutations with lucky set $[k]$. 

\begin{block}
Let $I=\{1,2,\ldots,k\}\subseteq [n]$. For each permutation $\pi=\pi_1\pi_2\cdots \pi_n \in \out_n(I;\, 213)$, there exists a permutation
\[\sigma=\sigma_1\sigma_2\cdots\sigma_k \in \Av_k(213),\]
whose first descent occurs at position $j$ ($\sigma_j>\sigma_{j+1}$), taking $j=k$ if $\sigma$ has no descents, such that $\pi$ admits a block decomposition
\[\pi=B_1B_2\cdots B_{k} \in \out_n(I;\, 213),\]
satisfying: 
\begin{enumerate}
    \item For any $r\in[k]$, the first element of the block $B_r$ is the lucky car $\sigma_r$.
    \item If $r>j$ then $B_r=\sigma_r$.
    \item The subsequence $\sigma_1\sigma_2\cdots\sigma_j$ is an increasing subsequence since $j$ is the position where the first descent occurs.
    \item For $r\leq j$, each block $B_r$ consists of the lucky car $\sigma_r$ followed by an increasing sequence (possibly empty) of unlucky cars from $[n]\setminus I=\{k+1,k+2,\ldots,n\}$.
    \item  For any $u,v\in [n]\setminus [k]$, if $u$ and $v$ are unlucky cars belonging to blocks $B_r$ and $B_s$, respectively, with $r<s\leq k$, then $u>v$. 
\end{enumerate}
\end{block}

\begin{example}
Let $n=9$ and $k=5$. The following are examples of parking permutations in $\Av_9(213)$ with lucky set $[5]$. The blocks are indicated with underbraces and each element at the beginning of each block is a lucky car.

\[
   \pi= \underbrace{189}_{B_1} \underbrace{2}_{B_2} \underbrace{567}_{B_3} \underbrace{3}_{B_4} \underbrace{4}_{B_5}
\]

\[ \sigma = \underbrace{49}_{B_1} \underbrace{5678}_{B_2} \underbrace{1}_{B_3} \underbrace{2}_{B_4} \underbrace{3}_{B_5}\]

\[ \tau = \underbrace{2}_{B_1} \underbrace{1}_{B_2} \underbrace{3}_{B_3} \underbrace{59}_{B_4} \underbrace{4678}_{B_5}\]
\end{example}

To verify the above decomposition, we begin with a lemma that we use for the characterization.

\begin{lemma}
\label{lem:descentsin213avoiding}
In a $213$-avoiding parking permutation with lucky cars $\{1,2,\ldots, k\}$, let $\sigma = \sigma_1 \sigma_2 \cdots \sigma_k$ be the subsequence consisting of all lucky cars. If there is a descent at position $j$ in $\sigma$, then all unlucky cars must appear in positions to the left of $\sigma_{j+1}$.
\end{lemma}

\begin{proof}
   Suppose $\sigma$ has a descent at position $j$, so $\sigma_j > \sigma_{j+1}$. Suppose to the contrary that there is some unlucky car $\tau$ appearing to the right of $\sigma_{j+1}$. Then since $\sigma_j \leq k$ and $\tau>k$, we have that $\tau>\sigma_j>\sigma_{j+1}$ and they form the subsequence $\sigma_j, \sigma_{j+1}, \tau$, which is order-isomorphic to $213$ in the parking permutation. This is a contradiction as the permutation must avoid $213$.
\end{proof}

In order to proceed with the enumeration the $213$-avoiding parking permutations, we enumerate the permutations that avoid $213$ whose first descent occurs at position $j$. The proof of the following lemma is similar to the recursive proof provided by Connolly et al.~\cite{firstascent}. 

\begin{lemma}
\label{lem:213Catalan}
The number of $213$-avoiding permutations of $[n]$ at which the first descent occurs at position $j$ is given by the $j$-fold Catalan convolution $C_{n,j}$.
\end{lemma}
\begin{proof}
The key observation here is that any $213$-avoiding permutation of $[n]$ becomes a $213$-avoiding permutation of $[n-1]$ when the element $n$ is removed from it. Hence, we can view each $213$-avoiding permutation of $[n]$ as being constructed by taking a $213$-avoiding permutation of $[n-1]$ and inserting $n$ at certain positions.

Suppose $\sigma\in \Av_{n-1}(213)$ with first descent in position $j$, so  $\sigma_j>\sigma_{j+1}$. If $n$ is inserted before $\sigma_1$ then $n\sigma_1\sigma_2\cdots\sigma_n$ avoids $213$ and the first descent is now at position $j=1$.

If $n$ is inserted anywhere after $\sigma_1$ and before $\sigma_{j+1}$, then the location at which $n$ is inserted becomes the position of the new first descent, and so the first descent is at position $\ell$ where $2\leq \ell \leq j+1$. Furthermore, as there are no descents before $n$, and $n$ cannot be involved in any subsequent descents, the new permutation is still $213$-avoiding. Finally, if $n$ is inserted after the original first descent, then the resulting permutation contains the pattern $213$ since the relative order of the subsequence $\sigma_j,\sigma_{j+1},n$ is $213$ because $n>\sigma_j>\sigma_{j+1}$. 

If $\sigma=12\cdots(n-1) = \id$, then there are no descents and we define $j=n-1$. 
With this convention, the above also holds when $\sigma$ is the identity permutation. 

Thus, a $213$-avoiding permutation of length $n-1$ with first descent at position $j$ gives rise to exactly one $213$-avoiding permutation of length $n$ with first descent at position $\ell$ for each $ 1\leq \ell \leq j+1$. 

Let $A(n,j)$ be the number of $213$-avoiding permutations of length $n$ with first descent at position $j$, then $A(n,j)$ is equal to the number of $213$-avoiding permutations of length $n-1$ with first descent at position $j-1$ or later, as it is these permutations that give rise to them. Thus,
\[A(n,j)=\sum_{\ell=j-1}^n A(n-1,\ell),\]
with initial conditions $A(n,0)=0$ and $A(n,n)=1$ (the identity permutation). Since the recurrence $A(n,j)$ coincides with the recurrence proven by Connolly et al.~\cite{firstascent}, the result follows. 
\end{proof}

\begin{lemma}
\label{lem:213descendingblocks}
    In a $213$-avoiding parking permutation with lucky cars $\{1,2,\ldots, k\}$, let $B_1$ and $B_2$ be contiguous blocks of unlucky cars separated by at least one lucky car, such that $B_1$ appears to the left of $B_2$. Then all entries of $B_1$ are strictly greater than all entries of $B_2$.
\end{lemma}

\begin{proof}
Let $B_1$, $B_2$ be as described above and let $\sigma_i$ be a lucky car such that $\sigma_i$ appears after $B_1$ and before $B_2$.
Suppose for the sake of contradiction that there exists some $u\in B_1$ and some $v\in B_2$ with $u<v$. Since $u$ and $v$ are both unlucky, we have that $v>u>k\geq \sigma_i$. Then $u,\sigma_i,v$ is a subsequence of the parking permutation which is order-isomorphic to 213, a contradiction.
\end{proof}

\begin{theorem}\label{thm:213avoid}[Catalan's triangle (read by rows), OEIS \seqnum{A009766}]
The number of parking permutations of length $n$ that avoid the pattern $213$ with lucky set $I=\{1,2,\ldots, k\}$ is
\[|\out_n(I;\, 213)|
= \frac{n-k+2}{n+1}\binom{n+k-1}{k-1}.
\]
\end{theorem}
\begin{proof}

To construct a parking permutation with lucky cars $\{1,2,\ldots, k\}$, we first select a permutation $\sigma \in \Av_k(213)$; this determines the placement of the lucky cars. Then, we place the unlucky cars $\{k+1,k+2,\ldots,n\}$. 

If $\sigma\neq 12\cdots k = \id$, let $j$ be the position of the first descent of $\sigma$. Otherwise if $\sigma=\id$, let $j=k$. By Lemma \ref{lem:descentsin213avoiding}, all unlucky cars may only be placed to the right of any of the lucky cars $\sigma_1,\sigma_2,\ldots,\sigma_{j}$. We can distribute them into the ``gaps'' after each lucky car with any number of unlucky cars in each gap. However, we must do this in such a way that the parking permutation still avoids $213$.

For each $i\in[j]$, let $\alpha_i$ denote the number of unlucky cars placed directly after $\sigma_i$. Note that the $\alpha_i$ can be 0 and the tuple $(\alpha_1,\alpha_2,\ldots, \alpha_j)$ forms a weak composition of $n-k$.
From Lemma \ref{lem:213descendingblocks} it follows that the contiguous block of unlucky cars placed after $\sigma_1$ must consist of the $\alpha_1$ unlucky cars of highest index, the contiguous block of unlucky cars placed after $\sigma_2$ must consist of the $\alpha_2$ unlucky cars of next highest index, and so on. Moreover, by Lemma \ref{lem:luckydescents}, each contiguous block of unlucky cars must be strictly increasing. So once the values of the $\alpha_i$ are chosen, the placement of the unlucky cars is completely determined. The number of such choices for the weak composition $(\alpha_1,\alpha_2,\ldots,\alpha_j)$ is given by $\binom{n-k+j-1}{j-1}$.

This count depends on the position of the first descent $j$ of the chosen $\sigma \in \Av_k(213)$ (if such a descent exists). This parameter $j$ can range from $1$ to $k-1$, and by Lemma \ref{lem:213Catalan}, the number of $\sigma \in \Av_k(213)$ with first descent at $j$ is given by $C_{k,j}$. Then in total, the number of $213$-avoiding parking permutations with lucky set $\{1,2,\ldots,k\}$ is
\[ \underbrace{\binom{n-k+k-1}{k-1}}_{\text{when $\sigma=\id$}}
+\underbrace{\sum_{j=1}^{k-1} \binom{n-k+j-1}{j-1}C_{k,j}}_{\text{when $\sigma$ has a descent}} = \sum_{j=1}^{k} \binom{n-k+j-1}{j-1}C_{k,j} = \frac{n-k+2}{n+1}\binom{n+k-1}{k-1},\]
where the first equality follows because $C_{k,k}=1$ and the second equality follows from Lemma~\ref{lem:213simplify}.
\end{proof}

\subsection{The class $\Av(312)$ with fixed lucky sets}
We begin with the block decomposition for the $312$-avoiding parking permutations with lucky set $[k]$. 

\begin{block}
Let $I=\{1,2,\ldots,k\}\subseteq [n]$. For each permutation $\pi=\pi_1\pi_2\cdots \pi_n \in \out_n(I;\, 312)$, there exists a permutation
\[\sigma=\sigma_1\sigma_2\cdots\sigma_k \in \Av_k(312),\]
whose last ascent occurs at position $j$ ($\sigma_j<\sigma_{j+1}$), taking $j=0$ if $\sigma$ has no ascents, such that $\pi$ admits a block decomposition
\[\pi=B_1B_2\cdots B_{k} \in \out_n(I; \, 312),\]
satisfying: 
\begin{enumerate}
    \item For any $r\in[k]$, the first element of the block $B_r$ is the lucky car $\sigma_r$.
    \item If $r<j$ then $B_r=\sigma_r$.
    \item The subsequence $\sigma_{j+1}\sigma_{j+2}\cdots\sigma_k$ is a decreasing subsequence since $j$ is the position where the last ascent occurs.
    \item For $r\geq j$, each block $B_r$ consists of the lucky car $\sigma_r$ followed by an increasing sequence (possibly empty) of unlucky cars from $[n]\setminus I=\{k+1,k+2,\ldots,n\}$.
    \item  For any $u,v\in [n]\setminus [k]$, if $u$ and $v$ are unlucky cars belonging to blocks $B_r$ and $B_s$, respectively, with $r<s\leq k$, then $u<v$. 
\end{enumerate}
\end{block}

\begin{example}
Let $n=9$ and $k=5$. The following are examples of parking permutations in $\Av_9(312)$ with lucky set $[5]$. The blocks are indicated with underbraces and each element at the beginning of each block is a lucky car.

\[
   \pi= \underbrace{26}_{B_1} \underbrace{5}_{B_2} \underbrace{478}_{B_3} \underbrace{3}_{B_4} \underbrace{19}_{B_5}
\]

\[ \sigma = \underbrace{3}_{B_1} \underbrace{2}_{B_2} \underbrace{16}_{B_3} \underbrace{5}_{B_4} \underbrace{4789}_{B_5}\]

\[ \tau = \underbrace{4}_{B_1} \underbrace{3}_{B_2} \underbrace{567}_{B_3} \underbrace{289}_{B_4} \underbrace{1}_{B_5}\]
\end{example}

\begin{lemma}[Correspondence between $213$- and $312$-avoiding permutations]\label{lemma:reverse bijection2}
Let $\pi=\pi_1\pi_2\cdots\pi_k\in \Av_k(213)$ whose first descent appears at position $j$ ($\pi_j>\pi_{j+1}$). Then, $\pi^r=\pi_k\pi_{k-1}\cdots\pi_1\in \Av_{k}(312)$ and the last ascent of $\pi^r$ is at position $k-j$ ($\pi_{k-j}<\pi_{k-j+1}$).
\end{lemma}

\begin{proof}
The proof is similar to Lemma~\ref{lemma:reverse bijection}.
\end{proof}

\begin{corollary}(Catalan's triangle (read by rows), OEIS \seqnum{A009766}) 
\label{cor:enumerationof312avoidingpermutationsbylastascent}
The number of $\sigma \in \Av_k(312)$ whose last ascent appears at position $j$ is given by $C_{k,k-j} = \frac{k-j}{k+j}\binom{k+j}{k}$.
\end{corollary}

\begin{proof}
The proof is similar to Corollary~\ref{cor:enumerationof321avoidingpermutationsbylastdescent}.
\end{proof}

\begin{lemma}
\label{lem:ascentsin312avoiding}
In a $312$-avoiding parking permutation with lucky cars $\{1,2,\ldots, k\}$, let $\sigma = \sigma_1 \sigma_2 \cdots \sigma_k$ be the subsequence consisting of all lucky cars. If there is an ascent at position $j$ in $\sigma$, then all unlucky cars must appear in positions to the right of $\sigma_{j}$.
\end{lemma}

\begin{proof}
Suppose $\sigma$ has an ascent at position $j$, so $\sigma_j > \sigma_{j+1}$. Suppose to the contrary that there is some unlucky car $\tau$ which appears to the left of $\sigma_{j}$. Then since $\sigma_j<\sigma_{j+1} \leq k$ and $\tau>k$, we have that $\tau>\sigma_{j+1}>\sigma_{j}$ and they form the subsequence $\tau, \sigma_j, \sigma_{j+1}$ which is order-isomorphic to $312$, which is a contradiction as the permutation must avoid $312$.
\end{proof}

\begin{lemma}
\label{lem:312ascendingblocks}
In a $312$-avoiding parking permutation with lucky cars $\{1,2,\ldots, k\}$, let $B_1$ and $B_2$ be contiguous blocks of unlucky cars separated by at least one lucky car, such that $B_1$ appears to the left of $B_2$. Then all entries of $B_1$ are strictly less than all entries of $B_2$.
\end{lemma}
\begin{proof}
Let $B_1,B_2$ be as described above and let $\sigma_i$ be a lucky car such that $\sigma_i$ appears after $B_1$ and before $B_2$. Suppose for the sake of contradiction that there exists some $u\in B_1$ and some $v\in B_2$ with $u>v$. Since $u$ and $v$ are both unlucky, we have that $u>v>k\geq \sigma_1$. Then $u,\sigma_i,v$ is a subsequence of the parking permutation, which is order-isomorphic to $312$, 
a contradiction.
\end{proof}

\begin{theorem}
\label{thm:312avoid}
The number of parking permutations of length $n$ that avoid the pattern $312$ with lucky set $I=\{1,2,\ldots, k\}$ is
\[|\out_n(I; \, 312)| = \frac{n-k+1}{n+1}\binom{n+k}{k} -\binom{n-1}{k}.\]
\end{theorem}

\begin{proof}

We follow a similar argument to the proof of Theorem~\ref{thm:213avoid}. 

We begin by selecting a permutation $\sigma \in \Av_k(312)$, which determines the placement of the lucky cars $\{1,2,\ldots, k\}$. We then choose how to place the remaining $n-k$ unlucky cars.

If $\sigma = k(k-1)\cdots 21=\id^r$, then the unlucky cars may be placed into any of the ``gaps'' directly after any of the lucky cars. (We cannot place them before the first lucky car, since the first car in the parking permutation is always lucky.)
In this case, we have $k$ gaps, and $n-k$ cars to distribute over them. We choose a weak composition $\alpha = (\alpha_1,\alpha_2,\ldots, \alpha_{k})$ of $n-k$ into $k$ parts, so we have $\binom{n-k+k-1}{k-1} =\binom{n-1}{k-1}$ choices for the number of unlucky cars to place into each of the gaps. 

Otherwise suppose $\sigma\neq \id^r$. Then let $j$ be the position of the last ascent of $\sigma$. 
By Lemma \ref{lem:ascentsin312avoiding}, all unlucky cars may only be placed to the right of any of the lucky cars $\sigma_j,\sigma_{j+1},\ldots,\sigma_k$. Again, we may place them into the ``gaps'' directly after each of these lucky cars, with any number of unlucky cars in each gap. 
For each $i\in[k-j+1]$, let $\alpha_i$ denote the number of unlucky cars placed directly after $\sigma_{j+i-1}$. Then the tuple $(\alpha_1,\alpha_2,\ldots, \alpha_{k-j+1})$ forms a weak composition of $n-k$, so there are $\binom{n-k+k-j+1-1}{k-j+1-1} = \binom{n-j}{k-j}$ choices for the number of unlucky cars to place into each of the gaps, when the last ascent of $\sigma$ occurs at position $j$.

Once we have determined how many unlucky cars to place in each gap, we must choose their order in such a way that the parking permutation still avoids $312$.
Lemma \ref{lem:312ascendingblocks} implies that the contiguous block of unlucky cars placed in the first gap from the left consists of the lowest $\alpha_1$ unlucky cars, the contiguous block of unlucky cars placed in the next gap consists of the next lowest $\alpha_2$ unlucky cars, and so on. Moreover, by Lemma \ref{lem:luckydescents}, each contiguous block of unlucky cars must be strictly increasing. So once the values of the $\alpha_i$ are chosen, the placement and ordering of the unlucky cars is completely determined.

The value of $j$, which records the position of the last ascent of the chosen $\sigma \in \Av_k(312)$ (if such an ascent exists), can range from $1$ to $k-1$, and by Corollary \ref{cor:enumerationof312avoidingpermutationsbylastascent}, the number of $\sigma \in \Av_k(312)$ with last ascent at position $j$ is given by $C_{k,k-j}$. Then the number of $312$-avoiding parking permutations with lucky set $\{1,2,\ldots,k\}$ is

\begin{align*}
|\out_n(I; \, 312)| = \underbrace{\binom{n-1}{k-1}}_{\text{when $\sigma=\id^r$}}
+\underbrace{\sum_{j=1}^{k-1} \binom{n-j}{k-j}C_{k,k-j}}_{\text{when $\sigma$ has an ascent}}&=\binom{n-1}{k-1}+\sum_{j=1}^{k-1}  \frac{k-j}{k+j}\binom{k+j}{k}\binom{n-j}{k-j}\\
&=\binom{n-1}{k-1}+\sum_{j=1}^{k-1}\frac{k-j}{k+j}\frac{k+j}{k}\frac{(k+j-1)!}{(k-1)!j!}\binom{n-j}{k-j}\\
&=\binom{n-1}{k-1}+\sum_{j=1}^{k-1} \frac{k-j}{k}\binom{k+j-1}{k-1}\binom{n-j}{k-j}.
\end{align*}
Now, by Lemma~\ref{lem:sumid1},
\[\sum_{j=1}^{k-1} \frac{k-j}{k}\binom{k+j-1}{k-1}\binom{n-j}{k-j}=\frac{n-k+1}{n+1}\binom{n+k}{k}-\binom{n}{k}.\]
Hence, 
\begin{align*}
|\out_n(I; \, 312)| &=\binom{n-1}{k-1}+ \frac{n-k+1}{n+1}\binom{n+k}{k}-\binom{n}{k}\\
&= \binom{n-1}{k-1} -\binom{n}{k} + \frac{n-k+1}{n+1}\binom{n+k}{k}\\
&= -\binom{n-1}{k} + \frac{n-k+1}{n+1}\binom{n+k}{k} &  \text{by Pascal's rule}\\
&= \frac{n-k+1}{n+1}\binom{n+k}{k} -\binom{n-1}{k}. & \qedhere
\end{align*}
\end{proof}

\subsection{The class $\Av(231)$ with fixed lucky sets}
We begin with the block decomposition for the $231$-avoiding parking permutations with lucky set $[k]$.

\begin{block}
Let $I=\{1,2,\ldots,k\}\subseteq [n]$. Each permutation $\pi=\pi_1\pi_2\cdots \pi_n \in \out_n(I;\, 231)$ admits a block decomposition
\[\pi=B_1B_2\cdots B_{\ell}\]
satisfying: 
\begin{enumerate}
    \item Each block $B_r$ avoids the pattern $231$.
    \item $B_1$ consists of only lucky cars.
    \item $B_\ell$ consists of a (possibly empty) increasing run of unlucky cars.
    \item For any $r\in\{2,3,\ldots,\ell-1\}$, either $B_r$ consists of exactly one unlucky car, or consists of only lucky cars.
     \item  For any $u,v\in [k]$, if $u$ and $v$ are lucky cars belonging to blocks $B_r$ and $B_s$, respectively, with $r<s$, then $u<v$. 
    \item  For any $u,v\in [n]\setminus[k]$, if $u$ and $v$ are unlucky cars belonging to blocks $B_r$ and $B_s$, respectively, with $r<s<\ell$, then $u>v$.
\end{enumerate}
\end{block}

\begin{example}
Let $n=9$ and let $k=5$. The following are examples of parking permutations in $\Av_9(231)$ with lucky set $[5]$. The blocks are indicated with underbraces.

\[
\pi= \underbrace{21}_{B_1} \underbrace{7}_{B_2} \underbrace{3}_{B_3} \underbrace{6}_{B_4} \underbrace{54}_{B_5} \underbrace{89}_{B_6}
\]

\[ \sigma = \underbrace{2134}_{B_1} \underbrace{8}_{B_2} \underbrace{5}_{B_3} \underbrace{679}_{B_4} \]

\[ \tau = \underbrace{1}_{B_1} \underbrace{9}_{B_2} \underbrace{2}_{B_3} \underbrace{8}_{B_4} \underbrace{3}_{B_5}
\underbrace{7}_{B_6}
\underbrace{4}_{B_7}
\underbrace{6}_{B_8}
\underbrace{5}_{B_9}
\underbrace{}_{B_{10}}\]

In $\tau$ above, $B_{10}$ is an empty block (there are no unlucky cars appearing after the last lucky car $5$).

\end{example}

\begin{lemma}
\label{lem:231unluckydecreasing}
     In a $231$-avoiding parking permutation with lucky cars $\{1,2,\ldots, k\}$, let $\sigma_k$ denote the rightmost lucky car. All unlucky cars appearing to the left of $\sigma_k$ must appear in decreasing order.
\end{lemma}

\begin{proof}
     Suppose to the contrary that unlucky cars $u,v \geq k$ both appear to the left of $\sigma_k$, with $u$ appearing to the left of $v$ and $u<v$. Since $\sigma_k$ is lucky, we also have $\sigma_k\leq k <u<v$. Then the subsequence $u,v,\sigma_k$ is order-isomorphic to 231, which is a contradiction.
\end{proof}

\begin{corollary}
\label{cor:231avoid}
    In a $231$-avoiding parking permutation with lucky cars $\{1,2,\ldots, k\}$, let $\sigma_k$ denote the rightmost lucky car. There cannot be two consecutive unlucky cars appearing to the left of $\sigma_k$.
\end{corollary}

\begin{proof}
    Suppose to the contrary that unlucky cars $u,v \geq k$ appear consecutively in the parking permutation somewhere to the left of $\sigma_k$.
    Since $u$ and $v$ appear consecutively and are both unlucky, by Lemma~\ref{lem:luckydescents}, $u$ and $v$ must appear in increasing order, but this contradicts Lemma \ref{lem:231unluckydecreasing}.
\end{proof}

\begin{lemma}
\label{lem:231ascendingblocks}
In a $231$-avoiding parking permutation with lucky cars $\{1,2,\ldots, k\}$, let $B_1$ and $B_2$ be contiguous blocks of lucky cars separated by at least one unlucky car, such that $B_1$ appears to the left of $B_2$. Then all entries of $B_1$ are strictly less than all entries of $B_2$.
\end{lemma}
\begin{proof}
Let $B_1,B_2$ be as described above and let $u$ be an unlucky car such that $u$ appears after $B_1$ and before $B_2$. Suppose for the sake of contradiction that there exists some $s\in B_1$ and some $t\in B_2$ with $s>t$. Since $s$ and $t$ are both lucky, we have that $t<s\leq k <u$. Then $s,u,t$ is a subsequence of the parking permutation which is order-isomorphic to $231$, a contradiction.
\end{proof}

\begin{lemma}\label{lem:catalan_prod}
For a composition $\alpha=(\alpha_1,\ldots, \alpha_{m+1}) \vdash k$,
\[\sum_{\alpha = (\alpha_1,\ldots, \alpha_{m+1}) \vdash k} C_{\alpha_1}C_{\alpha_2}\cdots C_{\alpha_{m+1}} = \frac{m+1}{k}\binom{2k}{k-m-1}. \]
\end{lemma}
\begin{proof}
Define \[S=\sum_{\alpha = (\alpha_1,\ldots, \alpha_{m+1}) \vdash k} C_{\alpha_1}C_{\alpha_2}\cdots C_{\alpha_{m+1}}.\] 
Let $C(z) = \sum_{n\geq 0 } C_nz^n$ be the generating function for the Catalan numbers. 
The Catalan convolution formula (see for instance~\cite[Equation 19]{LarcombeFrench} and \cite{regev2011proofcatalansconvolutionformula}) states
\[
\sum_{i_1+i_2+\cdots+i_k=n+k} C_{i_1-1}C_{i_2-1}\cdots C_{i_k-1} = \frac{k}{2n+k}\binom{2n+k}{n},
\]
with $i_j\geq 1$ for all $j$.
Since the left-hand side of the previous equation is the coefficient of $z^n$ in the expansion of $(C(z))^k$, (see for example~\cite[Equations 5.70 and 7.70]{GKP}), we have
\begin{equation}\label{eq:coeff-for-weakcomp}
    [z^n] C(z)^k =\frac{k}{2n+k}\binom{2n+k}{n}.
\end{equation}

Observe that the summation $S$ is over all compositions $\alpha=(\alpha_1,\ldots, \alpha_{m+1})\vdash k$ where each $\alpha_j \geq 1$. Hence, we compute the coefficient of $z^k$ in $(C(z)-C_0)^{m+1}$. 
The Catalan generating function satisfies $C(z)=1+zC(z)^2$. Thus,
\begin{align*}
[z^k]\left(\sum_{n\geq 1} C_nz^n\right)^{m+1} = [z^k]\left(\sum_{n\geq 0} C_nz^n- C_0\right)^{m+1}= [z^k](C(z)-1)^{m+1} &= [z^k]\left(zC(z)^2 \right)^{m+1}\\
&=[z^k]z^{m+1}C(z)^{2m+2}\\
&=[z^{k-m-1}]C(z)^{2m+2}.
\end{align*}
Now, making the substitutions $n\coloneqq k-m-1$ and $k\coloneqq 2m+2$ into Equation~\eqref{eq:coeff-for-weakcomp}, we obtain 
\[[z^k]\left(\sum_{n\geq 1} C_nz^n\right)^{m+1}  = \frac{m+1}{k}\binom{2k}{k-m-1}. \qedhere\]
\end{proof}

\begin{theorem}\label{thm:231avoid}[Catalan's triangle (read by rows), OEIS \seqnum{A009766}]
The number of parking permutations of length $n$ that avoid the pattern $231$ with lucky set $I=\{1,2,\ldots, k\}$ is
\[|\out_n(I;\, 231)|=\frac{n-k+2}{n+1}\binom{n+k-1}{k-1}. \]
\end{theorem}
\begin{proof}
Let $I=\{1,2,\ldots, k\}$ be a lucky set. 
We construct $231$-avoiding parking permutations first by placing the unlucky cars, and then placing the lucky cars in such a way that the resulting permutation still avoids $231$.

By Lemma~\ref{lem:231unluckydecreasing}, there exists a (possibly empty) subsequence of unlucky cars $\{k+1,k+2,\ldots, n\}$ that must appear in decreasing order, to the left of the rightmost lucky car $\sigma_k$. By Lemma~\ref{lem:luckydescents}, the remaining $n-k-m$ unlucky cars appearing after $\sigma_k$ will appear in increasing order.
Let $A$ be the subset of all unlucky cars to the left of $\sigma_k$, where $A$ contains $m$ cars, $0\leq m\leq n-k$. 
There are then $\binom{n-k}{m}$ ways to select the set $A$.

Now we place the lucky cars.
Since the cars in $A$ are those which are to appear before the rightmost lucky car, we must place at least one lucky car directly after the last car in the set $A$. In addition, since all of the cars in $A$ are unlucky, we must place at least one lucky car at the beginning of the parking permutation, directly before the first car in $A$. Furthermore, 
by Corollary~\ref{cor:231avoid},
at least one lucky car must be placed directly after each of the first $m-1$ cars in $A$.
Therefore, if we select $m$ unlucky cars to form the subset $A$, then there are $m+1$ locations in which we must insert at least one lucky car. Determining the number of lucky cars to place in each possible gap amounts to choosing a composition $\alpha=(\alpha_1,\alpha_2,\ldots, \alpha_{m+1})$ of $k$ into $m+1$ parts.

Once we have determined how many lucky cars to place in each gap, we must then order them in such a way that the parking permutation still avoids $231$.
From Lemma \ref{lem:231ascendingblocks}, it follows that the contiguous block of lucky cars placed in the first gap from the left consists of the lowest $\alpha_1$ lucky cars, the contiguous block of lucky cars placed in the next gap consists of the next lowest $\alpha_2$ lucky cars, and so on. So the set of lucky cars appearing in each gap is determined once the composition $\alpha$ is selected. The only restriction on the order of the lucky cars within each gap is that they must still respect the $231$ avoidance condition. So for each $i\in [m+1]$, the $\alpha_i$ cars in the $i$th gap may be placed in $C_{\alpha_i}$ possible ways. The number of choices for how to order all of the lucky cars is then $C_{\alpha_1}C_{\alpha_2}\cdots C_{\alpha_{m+1}}$.

Then the number of $231$-avoiding parking permutations with lucky set $\{1,2,\ldots,k\}$ is
\begin{align*}
    |\out_n(I;\, 231)| &=
\sum_{m=0}^{n-k} \binom{n-k}{m} \sum_{\alpha = (\alpha_1,\ldots, \alpha_{m+1}) \vdash k} C_{\alpha_1}C_{\alpha_2}\cdots C_{\alpha_{m+1}}\\
&=\sum_{m=0}^{n-k} \frac{m+1}{k}\binom{2k}{k-m-1}\binom{n-k}{m} & \text{by Lemma~\ref{lem:catalan_prod}} \\
&=\frac{n-k+2}{n+1}\binom{n+k-1}{k-1} & \text{by Lemma~\ref{lem:identity4}}. \qquad \qedhere
\end{align*}
\end{proof}

\subsection{The class $\Av(132)$ with fixed lucky sets}

\begin{block}
Let $I=\{1,2,\ldots,k\}\subseteq [n]$. Each permutation $\pi=\pi_1\pi_2\cdots \pi_n \in \out_n(I;\, 132)$ admits a 
block decomposition
\[\pi=B_1B_2\cdots B_{\ell}\]
satisfying: 
\begin{enumerate}
    \item Each block $B_r$ avoids the pattern $132$.
    \item $B_1$ consists of only lucky cars with consecutive indices.
    \item For any $r\in\{2,3,\ldots,\ell\}$, either $B_r$ consists of exactly one unlucky car, or consists of only lucky cars with consecutive indices.
    \item  For any $u,v\in [k]$, if $u$ and $v$ are lucky cars belonging to blocks $B_r$ and $B_s$, respectively, with $r<s$, then $u>v$.
    \item For any $u,v\in [n]\setminus [k]$, if $u$ and $v$ are unlucky cars belonging to blocks $B_r$ and $B_s$, respectively, with $r<s\leq \ell$, then $u<v$.
\end{enumerate}
\end{block}

\begin{example}
Let $n=9$ and let $k=5$. The following are examples of parking permutations in $\Av_9(132)$ with lucky set $[5]$. The blocks are indicated with underbraces.

\[
   \pi= \underbrace{45}_{B_1} \underbrace{6}_{B_2} \underbrace{7}_{B_3} \underbrace{321}_{B_4} \underbrace{8}_{B_5} \underbrace{9}_{B_6}
\]

\[ \sigma = \underbrace{5}_{B_1} \underbrace{6}_{B_2} \underbrace{34}_{B_3} \underbrace{7}_{B_4} \underbrace{8}_{B_5} \underbrace{9}_{B_6} \underbrace{21}_{B_7} \]

\[ \tau = \underbrace{12345}_{B_1} \underbrace{6}_{B_2} \underbrace{7}_{B_3} \underbrace{8}_{B_4} \underbrace{9}_{B_5}\]

\end{example}

\begin{lemma}\label{lem:132increasing}
In a $132$-avoiding parking permutation with lucky set $\{1,2,\ldots,k\}$, the unlucky cars all appear in increasing order. 
\end{lemma}
\begin{proof}
Assume for contradiction that an unlucky car $s$ appears to the left of an unlucky car $t$ where $s>t$. Since $s$ is unlucky, it appears after some lucky car $v$. Note that $v\leq k$ and $s>t>k$. Thus $v<t<s$ and $v$ appears to the left of $s$ which appears to the left of $t$, forming a subsequence which is order-isomorphic to $132$. This is a contradiction.
\end{proof}

\begin{lemma}
\label{lem:132descendingblocks}
    In a $132$-avoiding parking permutation with lucky cars $\{1,2,\ldots, k\}$, let $B_1$ and $B_2$ be contiguous blocks of lucky cars which are separated by at least one unlucky car, with $B_1$ appearing to the left of $B_2$. Then all entries of $B_1$ are strictly greater than all entries of $B_2$.
\end{lemma}

\begin{proof}
Let $B_1$, $B_2$ be as described, and let $\tau$ be an unlucky car such that each car in $B_1$ appears to the left of $\tau$ and each car in $B_2$ appears to the right of $\tau$.
Suppose for contradiction that there exists some $u\in B_1$ and some $v\in B_2$ with $u<v$. Since $u$ and $v$ are both lucky, we have that $u<v\leq k<\tau$. Then $u,\tau,v$ form a subsequence of the parking permutation which is order-isomorphic to $132$, which is a contradiction.
\end{proof}

\begin{theorem}\label{thm:132avoid}[Catalan's triangle (read by rows), OEIS \seqnum{A009766}]
The number of parking permutations of length $n$ that avoid the pattern $132$ with lucky set $I=\{1,2,\ldots,k\}$ is
\[|\out_n(I;\, 132)|=\frac{n-k+2}{n+1}\binom{n+k-1}{k-1}. \]
\end{theorem}

\begin{proof}
Let $I=\{1,2,\ldots,k\}$ be a lucky set.
    We construct $132$-avoiding parking permutations first by placing the unlucky cars, and then placing the lucky cars in such a way that the resulting permutation still avoids $132$.

    By Lemma \ref{lem:132increasing}, the unlucky cars $\{k+1,k+2,\ldots, n\}$ must all appear in increasing order. So there is no choice regarding the ordering of the subsequence of unlucky cars.

    Next, the lucky cars may be placed in the ``gaps'' to the left of any of the unlucky cars, or after the last unlucky car. There are $n-k+1$ such gaps, and suppose that $\ell$ of them will have at least one lucky car added. Here, $\ell$ is a positive integer which can range from $1$ to $n-k+1$. Since the first car in the parking permutation is always lucky, we must place at least one lucky car in the leftmost gap. Then among the remaining $n-k$ gaps, we choose $\ell-1$ to place lucky cars into. There are $\binom{n-k}{\ell-1}$ choices for this step.

    Once the locations of the lucky cars are selected, we now must decide how many lucky cars to place into each gap. This amounts to choosing a composition $\alpha = (\alpha_1,\alpha_2,\ldots, \alpha_\ell)$ of $k$ into $\ell$ parts, where $\alpha_i$ determines the number of lucky cars deposited into gap $i$. 
    {From Lemma \ref{lem:132descendingblocks}, we have that the contiguous block of lucky cars placed in the first chosen gap contains the highest $\alpha_1$ lucky cars, the contiguous block of unlucky cars placed in the next chosen gap contains the next highest $\alpha_2$ lucky cars, and so on up until the $\ell$th chosen gap.} 

    So the set of lucky cars appearing in each gap is determined by the choice of the composition $\alpha$. Within each gap, the lucky cars can be ordered in any way that avoids 132. So for each $i\in [\ell]$, the $\alpha_i$ cars in the $i$th gap may be placed in $C_{\alpha_i}$ possible ways. Then the number of ways to order all of the lucky cars is $C_{\alpha_1}C_{\alpha_2}\cdots C_{\alpha_\ell}$.
    In total, we have 
    \begin{align*}
|\out_n(I;\, 132)|&=\sum_{\ell = 1}^{n-k+1} \binom{n-k}{\ell-1} \sum_{(\alpha_1,\alpha_2,\ldots, \alpha_\ell) \vdash k} C_{\alpha_1}  C_{\alpha_2}\cdots  C_{\alpha_\ell}\\
&= \sum_{\ell = 0}^{n-k} \binom{n-k}{\ell} \sum_{(\alpha_1,\alpha_2,\ldots, \alpha_{\ell+1}) \vdash k} C_{\alpha_1}  C_{\alpha_2}\cdots  C_{\alpha_{\ell+1}}\\
&= \sum_{\ell = 0}^{n-k} \frac{\ell +1}{k}\binom{2k}{k-\ell -1}\binom{n-k}{\ell} &\text{by Lemma~\ref{lem:catalan_prod}}\\
&= \frac{n-k+2}{n+1}\binom{n+k-1}{k-1} & \text{by Lemma~\ref{lem:identity4}}. \qquad \qedhere
\end{align*}
\end{proof}

\section{Catalan Convolution and Binomial Identities}

In some of our results, we use the following generalization of Vandermonde's identity; see~\cite{Gould}. 
\begin{proposition}{(Generalization of Vandermonde's convolution~\cite[Equation 3]{Gould})}\label{prop:GenVandermondes}

For integers $n,\alpha,\gamma$,
\[\sum_{k=0}^n \binom{\alpha + k}{k}\binom{\gamma + n-k}{n-k}=\binom{\alpha+\gamma +n+1}{n}.\]    
\end{proposition}

\begin{lemma}\label{lem:identity1}
For integers $n,k\geq 1$,
\[\sum_{i=n-k}^k \binom{i}{n-k}C_{k,i}=\binom{2k}{2k-n}-\binom{2k}{2k-n-1} = \frac{2n-2k+1}{n+1} \binom{2k}{n}.\]
\end{lemma}
\begin{proof}
We use the fact that $C_{k,i} = \frac{i}{2k-i} \binom{2k-i}{k}$. Then
\begin{align*}
    \sum_{i=n-k}^k \binom{i}{n-k}C_{k,i} &= \sum_{i=n-k}^k \binom{i}{n-k} \frac{i}{2k-i} \binom{2k-i}{k}.
\end{align*}
Rewrite $\frac{i}{2k-i}=\frac{2k}{2k-i}-1$ so that
\begin{align*}
  \sum_{i=n-k}^k \binom{i}{n-k} \frac{i}{2k-i} \binom{2k-i}{k} &= \sum_{i=n-k}^k \frac{2k}{2k-i} \binom{i}{n-k} \binom{2k-i}{k} - \sum_{i=n-k}^k \binom{i}{n-k} \binom{2k-i}{k}.
\end{align*}
Define
\begin{equation*}
    S(n,k)=\sum_{i=n-k}^k \frac{2k}{2k-i} \binom{i}{n-k} \binom{2k-i}{k}, \quad \text{ and } \quad 
T(n,k)=\sum_{i=n-k}^k \binom{i}{n-k} \binom{2k-i}{k}.
\end{equation*}
It suffices to show that $S(n,k)-T(n,k)$ is the claimed formula.

We have,
\begin{align*}
S(n,k)&= 2k \sum_{i=n-k}^k \frac{1}{2k-i} \binom{i}{n-k} \frac{(2k-i)!}{k!(k-i)!}\\&= 2\sum_{i=n-k}^k \binom{i}{n-k} \frac{(2k-i-1)!}{(k-1)!(k-i)!}\\
&=2\sum_{i=n-k}^k \binom{i}{n-k} \binom{2k-i-1}{k-1}\\
&=2\sum_{i=0}^{2k-n} \binom{n-k+i}{n-k}\binom{2k-(n-k+i)-1}{k-1} & \text{by reindexing}\\
&= 2\sum_{i=0}^{2k-n} \binom{n-k+i}{i}\binom{(k-1)+(2k-n)-i}{2k-n-i} &  \text{by symmetry of binomial coefficients}\\
&=2\binom{(n-k)+(k-1)+(2k-n)+1}{2k-n} & \text{by Vandermonde's identity}\\
&=2\binom{2k}{2k-n},
\end{align*}
where in the second-to-last step, we use the generalization of Vandermonde's identity from Proposition~\ref{prop:GenVandermondes} with the substitutions $n\coloneqq 2k-n, \alpha\coloneqq n-k, \gamma\coloneqq k-1$ and $k\coloneqq i $.
Next, we find a closed formula for $T(n,k)$.
Observe that
\begin{align*}
T(n,k)& = \sum_{i=n-k}^k \binom{i}{n-k} \binom{2k-i}{k}\\
&= \sum_{i=0}^{2k-n} \binom{n-k+i}{n-k}\binom{2k-(n-k+i)}{k} & \text{by reindexing}\\
&= \sum_{i=0}^{2k-n} \binom{n-k+i}{n-k}\binom{k+(2k-n)-i}{k}\\
&= \sum_{i=0}^{2k-n} \binom{n-k+i}{i}\binom{k+(2k-n)-i}{(2k-n)-i} & \text{by symmetry of binomial coefficients}\\
&=\binom{n-k+k+2k-n+1}{2k-n}  & \text{by Vandermonde's identity}\\
&=\binom{2k+1}{2k-n},
\end{align*}
where in the second-to-last step, we again apply the generalization of Vandermonde's identity from Proposition~\ref{prop:GenVandermondes}, here with the substitutions $n\coloneqq 2k-n, \alpha\coloneqq n-k, \gamma\coloneqq k$ and $k\coloneqq i $.

\noindent Altogether, we obtain
\begin{align*}
    S(n,k)-T(n,k)&= 2\binom{2k}{2k-n}- \binom{2k+1}{2k-n}\\
    &= \binom{2k}{2k-n} + \binom{2k}{2k-n} -\binom{2k+1}{2k-n}\\
    &=\binom{2k}{2k-n}-\binom{2k}{2k-n-1} & \text{by Pascal's rule}\\
       & = \frac{(2k!)}{(2k-n)!n!} \cdot \frac{(n+1)}{(n+1)}- \frac{(2k!)}{(2k-n-1)!(n+1)!}\cdot \frac{(2k-n)}{(2k-n)} \\
    &= \frac{(2k)!(n+1) - (2k)!(2k-n)}{(2k-n)!(n+1)!}\\
    &= \frac{(2k)! (2n-2k+1)}{(2k-n)!(n+1)n!}= \binom{2k}{n}\frac{2n-2k+1}{n+1}. &  \qedhere
\end{align*}
\end{proof}

\begin{lemma}
\label{lem:213simplify}
For positive integers $k,n$,
\[\sum_{j=1}^{k} \binom{n-k+j-1}{j-1}C_{k,j}= \frac{n-k+2}{n+1}\binom{n+k-1}{k-1}.\]
\end{lemma}
\begin{proof}

The argument for this identity follows a similar structure to that of Lemma~\ref{lem:identity1}.
We first make the substitution $C_{k,j} = \frac{j}{2k-j}\binom{2k-j}{k}$.
Then
\begin{align*}
    \sum_{j=1}^{k} \binom{n-k+j-1}{j-1}C_{k,j} &= \sum_{j=1}^{k} \binom{n-k+j-1}{j-1}\frac{j}{2k-j}\binom{2k-j}{k}.
\end{align*}
Rewriting the fraction $\frac{j}{2k-j} = \frac{2k}{2k-j}-1$, this becomes
\begin{align*}
    \sum_{j=1}^{k} \binom{n-k+j-1}{j-1}\frac{2k}{2k-j}\binom{2k-j}{k} - \sum_{j=1}^{k} \binom{n-k+j-1}{j-1}\binom{2k-j}{k}.
\end{align*}
Let $S(n,k): = \sum_{j=1}^{k} \binom{n-k+j-1}{j-1}\frac{2k}{2k-j}\binom{2k-j}{k} $ and let $T(n,k) = \sum_{j=1}^{k} \binom{n-k+j-1}{j-1}\binom{2k-j}{k}$. We now find closed formulas for these expressions. First,
\begin{align*}
S(n,k)= \sum_{j=1}^{k} \binom{n-k+j-1}{j-1}\frac{2k}{2k-j}\binom{2k-j}{k} &= \sum_{j=1}^{k} \binom{n-k+j-1}{j-1}\frac{2k}{2k-j}\frac{(2k-j)!}{k!(k-j)!} \\
&=2 \sum_{j=1}^{k} \binom{n-k+j-1}{j-1}\frac{(2k-j-1)!}{(k-1)!(k-j)!} \\
&=2 \sum_{j=1}^{k} \binom{n-k+j-1}{j-1}\binom{2k-j-1}{k-1}.
\end{align*}
Now we reindex so that $j$ is initially set to $0$ and by symmetry of binomial coefficients:
\begin{align*}
2 \sum_{j=0}^{k-1} \binom{n-k+j}{j}\binom{2k-j-2}{k-1} = 2 \sum_{j=0}^{k-1} \binom{n-k+j}{j}\binom{2k-j-2}{k-j-1}.
\end{align*}
Then using the generalized Vandermonde identity (Proposition \ref{prop:GenVandermondes}) with the substitutions $k:=j$, $n:=k-1$, $\gamma:= k-1$, $\alpha:= n-k$, we have
\[S(n,k) = 2 \binom{n-k+k-1+k-1+1}{k-1} = 2\binom{n+k-1}{k-1}.\]
Now for
\begin{align*}
    T(n,k)&=\sum_{j=1}^{k} \binom{n-k+j-1}{j-1}\binom{2k-j}{k},
\end{align*}
we begin by reindexing so that $j$ is initialized to 0:
\begin{align*}
\sum_{j=0}^{k-1} \binom{n-k+j}{j}\binom{2k-j-1}{k} =\sum_{j=0}^{k-1} \binom{n-k+j}{j}\binom{2k-j-1}{k-j-1} \qquad \text{by symmetry of binomial coefficients.}
\end{align*}

We now use Proposition \ref{prop:GenVandermondes} with the substitutions $k:=j$, $n:=k-1$, $\gamma = k-1$, $\alpha:=n-k$, $\gamma:=k$. This yields
$
T(n,k) = \binom{n+k}{k-1}$.
So 
\begin{align*}
\sum_{j=1}^{k} \binom{n-k+j-1}{j-1}C_{k,j} &= S(n,k) - T(n,k) \\
&= 2\binom{n+k-1}{k-1} - \binom{n+k}{k-1}\\
&= \binom{n+k-1}{k-1}+\binom{n+k-1}{k-1} - \binom{n+k}{k-1}\\
&=  \binom{n+k-1}{k-1}-\binom{n+k-1}{k-2}&\text{by Pascal's rule}\\
&=  \frac{(n+k-1)!}{(n)!(k-1)!}-\frac{(n+k-1)!}{(n+1)!(k-2)!} \\
&=  \frac{(n+k-1)!}{(n)!(k-1)!} \cdot \frac{n+1}{n+1}-\frac{(n+k-1)!}{(n+1)!(k-2)!}\cdot\frac{k-1}{k-1}\\
&=  \frac{(n+k-1)! (n+1 - (k-1))}{(n+1)(n)!(k-1)!}\\
&= \frac{n-k+2}{n+1}\binom{n+k-1}{k-1}. &\qedhere
\end{align*}
\end{proof}

\begin{lemma} For integers $k,n\geq 1$ and $i\geq 0$,
\label{lem:321simplify}
\[\sum_{j=1}^{k-1}  \binom{n-i-j-1}{k-j-1}C_{k,k-j}=\frac{n -i +2-k }{n -i +1} \binom{n-i+k-1}{k-1}-\binom{n -i -1}{k -1}.\]
\end{lemma}

\begin{proof}
    We begin by reindexing the sum $\sum_{j=1}^{k-1}  \binom{n-i-j-1}{k-j-1}C_{k,k-j}$. Define $m:=k-j$. Then $m$ ranges from $k-1$ (when $j=1$) down to $1$ (when $j=k-1$). Then
    \[\sum_{j=1}^{k-1}  \binom{n-i-j-1}{k-j-1}C_{k,k-j} = \sum_{m=1}^{k-1}  \binom{n-i-(k-m)-1}{m-1}C_{k,m}\]

    Next, we compute a formula for \[\sum_{m=1}^{k}  \binom{n-i-(k-m)-1}{m-1}C_{k,m} = \sum_{m=1}^{k}  \binom{n-i-k+m-1}{m-1}C_{k,m}.\] 
    Let $N := n-i$. Then this sum can be rewritten as
    \[
    \sum_{m=1}^{k}  \binom{N-k+m-1}{m-1}C_{k,m}
    \]
    By Lemma~\ref{lem:213simplify}, 
    this simplifies to
    \[
    \frac{N-k+2}{N+1}\binom{N+k-1}{k-1} = \frac{n-i-k+2}{n-i+1}\binom{n-i+k-1}{k-1}.
    \]
    Then going back to the original summation, we have
    \begin{align*}
        &\sum_{j=1}^{k-1}  \binom{n-i-j-1}{k-j-1}C_{k,k-j}\\
        &= \frac{n-i-k+2}{n-i+1}\binom{n-i+k-1}{k-1} - \sum_{m=k}^{k}  \binom{n-i-k+m-1}{m-1}C_{k,m}\\
        &=  \frac{n-i-k+2}{n-i+1}\binom{n-i+k-1}{k-1}  - \binom{n-i-k+k-1}{k-1}C_{k,k} \\
        &=  \frac{n-i-k+2}{n-i+1}\binom{n-i+k-1}{k-1} - \binom{n-i-1}{k-1}C_{k,k} \\
        &=  \frac{n-i-k+2}{n-i+1}\binom{n-i+k-1}{k-1} - \binom{n-i-1}{k-1} &\text{since $C_{k,k}= 1$.} \qquad \qedhere
    \end{align*}
\end{proof}

\begin{lemma}\label{lem:identity4}
For positive integers $n,k$,
\[\sum_{m=0}^{n-k} \frac{m+1}{k}\binom{2k}{k-m-1}\binom{n-k}{m}=\frac{n-k+2}{n+1}\binom{n+k-1}{k-1}. 
\]
\end{lemma}
\begin{proof}
Note that whenever $k-m-1<0$, $\binom{2k}{k-m-1}=0$. Thus, it suffices to show that 
\begin{equation}\label{eq:sumidentity1}
    \sum_{m=0}^{k-1} \frac{m+1}{k}\binom{2k}{k-m-1}\binom{n-k}{m}=\frac{n-k+2}{n+1}\binom{n+k-1}{k-1}.
\end{equation}
Notice that
\[\frac{m+1}{k}\binom{2k}{k-m-1}=\binom{2k}{k-m-1}-\frac{k-m-1}{k}\binom{2k}{k-m-1}\]
since $m+1=k-(k-m-1)$. Using that, we can rewrite the left hand side of Equation~\eqref{eq:sumidentity1} into
\begin{equation}\label{eq:sumidentity2}
    \sum_{m=0}^{k-1} \binom{n-k}{m}\binom{2k}{k-m-1} - \sum_{m=0}^{k-1} \binom{n-k}{m}\frac{k-m-1}{k}\binom{2k}{k-m-1}.
\end{equation}
Furthermore, we can rewrite the following term
\[ \frac{k-m-1}{k}\binom{2k}{k-m-1}=2\binom{2k-1}{k-m-2}.\]
Therefore, Equation~\eqref{eq:sumidentity2} becomes
\begin{equation}\label{eq:sumidentity3}
    \sum_{m=0}^{k-1} \binom{n-k}{m}\binom{2k}{k-m-1} - 2\sum_{m=0}^{k-1} \binom{n-k}{m}\binom{2k-1}{k-m-2}.
\end{equation}
We use the generalization of Vandermonde's convolution (Proposition~\ref{prop:GenVandermondes}) to obtain
\[\sum_{m=0}^{k-1} \binom{n-k}{m}\binom{2k}{k-m-1} =\binom{n-k+2k}{k-1}=\binom{n+k}{k-1},\]
and \[ \sum_{m=0}^{k-1} \binom{n-k}{m}\binom{2k-1}{k-m-2}=\binom{n-k+2k-1}{k-2}=\binom{n+k-1}{k-2}.\]
Hence, Equation~\eqref{eq:sumidentity3} becomes
\begin{equation*}
\binom{n+k}{k-1}-2\binom{n+k-1}{k-2}=\frac{n+k}{n+1}\binom{n+k-1}{k-1}-\frac{2(k-1)}{n+1}\binom{n+k-1}{k-1}=\frac{n-k+2}{n+1}\binom{n+k-1}{k-1},
\end{equation*}
as desired.
\end{proof}

\begin{lemma}\label{lem:sumid1}
For any positive integers $k$ and $n$, \[\sum_{j=1}^{k-1} \frac{k-j}{k}\binom{k+j-1}{k-1}\binom{n-j}{k-j}=\frac{n-k+1}{n+1}\binom{n+k}{k}-\binom{n}{k}. \]
\end{lemma}
\begin{proof}
Observe that
\begin{align*}
 \frac{k-j}{k}\binom{k+j-1}{k-1}\binom{n-j}{k-j} &= \frac{1}{k}\frac{(k+j-1)!}{(k-1)!j!}(k-j)\frac{(n-j)!}{(k-j)!(n-k)!}\\
&= \frac{1}{j}\frac{(k+j-1)!}{k!(j-1)!}\frac{(n-j)!}{(k-j-1)!(n-k+1)!}(n-k+1)\\
&= \frac{n-k+1}{j}\binom{k+j-1}{k}\binom{n-j}{k-j-1}.
\end{align*}
Hence,
\begin{align*}
    \sum_{j=1}^{k-1} \frac{k-j}{k}\binom{k+j-1}{k-1}\binom{n-j}{k-j} &=(n-k+1)\sum_{j=1}^{k-1}\frac{1}{j}\binom{k+j-1}{k}\binom{n-j}{k-j-1}=(n-k+1)S,
\end{align*}
with
\begin{equation}\label{eq:Side}
S=\sum_{j=1}^{k-1}\frac{1}{j}\binom{k+j-1}{k}\binom{n-j}{k-j-1}.    
\end{equation}
It suffices to show that 
\begin{equation}\label{eq:iden-to-prove1}
S=\frac{1}{n+1}\binom{n+k}{k}-\frac{1}{n-k+1}\binom{n}{k}.
\end{equation}
Now, from Equation~\eqref{eq:Side},
\begin{align*}
\binom{k+j-1}{k}=\frac{k+j-1}{k}\frac{(k+j-2)!}{(k-1)!(j-1)!}=\frac{k+j-1}{k}\binom{k+j-2}{k-1},
\end{align*}
and
\begin{align*}
\binom{k+j-2}{k-1}=\frac{(k+j-2)!}{(k-1)!(j-1)!}=\frac{j}{k+j-1}\frac{(k+j-1)!}{(k-1)!j!}=\frac{j}{k+j-1}\binom{k+j-1}{k-1}.
\end{align*}
Altogether, 
\begin{align*}
 \binom{k+j-1}{k}=\frac{k+j-1}{k}\binom{k+j-2}{k-1}= \frac{k+j-1}{k}\frac{j}{k+j-1}\binom{k+j-1}{k-1}=\frac{j}{k}\binom{k+j-1}{k-1}.
\end{align*}
It follows that
\[
S=\sum_{j=1}^{k-1}\frac{1}{j}\binom{k+j-1}{k}\binom{n-j}{k-j-1}= \sum_{j=1}^{k-1}\frac{1}{j} \frac{j}{k}\binom{k+j-1}{k-1}\binom{n-j}{k-j-1}=\frac{1}{k} \sum_{j=1}^{k-1}\binom{k+j-1}{k-1}\binom{n-j}{k-j-1}.
\]
Equivalently,
\begin{equation}\label{eq:kSdef}
 kS= \sum_{j=1}^{k-1}\binom{k+j-1}{k-1}\binom{n-j}{k-j-1}.
\end{equation}
Next, we show that 
\begin{equation}\label{eq:idtoprove2}
\sum_{j=1}^{k-1}\binom{k+j-1}{k-1}\binom{n-j}{k-j-1}=\binom{n+k}{n+1}-\binom{n}{k-1}.
\end{equation}
Equivalently, since $\binom{k-1}{k-1}\binom{n}{k-1}=\binom{n}{k-1}$, then we proceed to show the following,
\begin{equation}\label{eq:sum2}
\sum_{j=0}^{k-1}\binom{k+j-1}{k-1}\binom{n-j}{k-j-1}=\binom{n+k}{n+1}.
\end{equation}
Let $\alpha\coloneqq k-1, \gamma\coloneqq n-k+1, n\coloneqq k-1$, and $k\coloneqq j$. Then, substituting those values into the generalized Vandermonde's convolution (Proposition~\ref{prop:GenVandermondes}), we obtain
\[
\sum_{k=0}^n \binom{\alpha + k}{k}\binom{\gamma + n-k}{n-k} = \sum_{j=0}^{k-1}\binom{k-1+j}{j}\binom{n-k+1+k-1-j}{k-1-j}= \binom{k-1+n-k+1+k-1+1}{k-1}.
\]
Simplifying, we get
\[\sum_{j=0}^{k-1}\binom{k-1+j}{j}\binom{n-j}{k-1-j}\\
= \binom{n+k}{k-1}=\binom{n+k}{n+1}.\]
Thus, we have proven that Equation~\eqref{eq:sum2} holds and so Equation~\eqref{eq:idtoprove2} follows. 

From Equation~\eqref{eq:kSdef}, we then obtain
\[ S= \frac{1}{k}\left(\binom{n+k}{n+1}-\binom{n}{k-1}\right).\]
Further simplification shows that 
\[S=\frac{1}{n+1}\binom{n+k}{k}-\frac{1}{n-k+1}\binom{n}{k},\]
which proves Equation~\eqref{eq:iden-to-prove1}, and hence the result follows.
\end{proof}

\section{Future work}

Many of the enumeration formulas in this work were found combinatorially as a summation involving Catalan numbers, and then algebraically simplified to more concise forms. It would be interesting to find a direct combinatorial explanation for the simplified enumeration formulas in terms of the parking permutations.

Another remaining question here is the enumeration of the pattern avoidance classes of parking permutations for lucky sets not of the form $\{1,2,\ldots, k\}$. That is, given an arbitrary lucky set $I$ and a pattern $\pi$, what is $|\out_n(I;\pi)|$? Alternatively, one could extend this study by keeping the restriction that the lucky set is of the form $[k]$ and computing $|\out_n([k], \pi)|$ for patterns $\pi$ of length greater than 3.

There are several generalizations and variations of parking permutations for which one could ask more pattern avoidance questions with refinement by the $\lucky$ statistic. For example, in an $(m,n)$-parking function, there can be more parking spots than cars, leaving some ``empty spaces'' in the corresponding parking permutation. One could define a notion of pattern avoidance in outcomes of $(m,n)$-parking functions and study the enumeration of such outcomes according to the set of lucky cars.

\section{Tool and computational resources disclosure}
The authors did not use AI or LLM tools to assist with any proofs or writing of this manuscript, with the exception of Overleaf language suggestions. Maple and SageMath were used for computations.

\section*{Acknowledgements}
L.~Martinez was supported by the NSF Graduate Research Fellowship Program under Grant No. 2233066.

\appendix
\section{Triangles of Integers}\label{appendix}

\begin{table}[H]
\begin{tabular}{lllllllllllllllll}
     &     &      &      &      &      &      & 1   &      & 1   &     &    &     &    &    &   &   \\
     &     &      &      &      &      & 2    &     & 3    &     & 1   &    &     &    &    &   &   \\
     &     &      &      &      & 5    &      & 9   &      & 5   &     & 1  &     &    &    &   &   \\
     &     &      &      & 14   &      & 28   &     & 20   &     & 7   &    & 1   &    &    &   &   \\
     &     &      & 42   &      & 90   &      & 75  &      & 35  &     & 9  &     & 1  &    &   &   \\
     &     & 132  &      & 297  &      & 275  &     & 154  &     & 54  &    & 11  &    & 1  &   &   \\
     & 429 &      & 1001 &      & 1001 &      & 637 &      & 273 &     & 77 &     & 13 &    & 1 &   \\
1430 &     & 3432 &      & 3640 &      & 2548 &     & 1260 &     & 440 &    & 104 &    & 15 &   & 1
\end{tabular}
\caption{(Lobb numbers, OEIS \seqnum{A039599}) The values $T(k,n)$ for $1\leq k\leq 9$ and $k\leq n\leq 2k$ such that $T(k,n)$ corresponds to the number of $123$-avoiding parking permutations of length $n$ with lucky set $I=[k]$. Here each row corresponds to a fixed value of $k$, and the entries in that row are ordered with $n$ ranging from $k$ to $2k$. }
\label{tab:123}
\end{table}

\begin{table}[H]
\centering
\begin{tabular}{lllllllllllllllll}
  &   &     &     &     &     &      &     & 1    &      &      &      &      &      &      &      &      \\
  &   &     &     &     &     &      & 1   &      & 2    &      &      &      &      &      &      &      \\
  &   &     &     &     &     & 1    &     & 4    &      & 5    &      &      &      &      &      &      \\
  &   &     &     &     & 1   &      & 8   &      & 13   &      & 14   &      &      &      &      &      \\
  &   &     &     & 1   &     & 16   &     & 32   &      & 41   &      & 42   &      &      &      &      \\
  &   &     & 1   &     & 32  &      & 76  &      & 112  &      & 131  &      & 132  &      &      &      \\
  &   & 1   &     & 64  &     & 176  &     & 292  &      & 380  &      & 428  &      & 429  &      &      \\
  & 1 &     & 128 &     & 400 &      & 736 &      & 1049 &      & 1290 &      & 1429 &      & 1430 &      \\
1 &   & 256 &     & 896 &     & 1808 &     & 2792 &      & 3698 &      & 4424 &      & 4861 &      & 4862
\end{tabular}
\caption{(New sequence) The values $T(n,k)$ for $1\le n\le 9$ and $1\le k\le n$ such that $T(n,k)$ corresponds to the number of $321$-avoiding parking permutations of length $n$ with lucky set $I=[k]$. Here each row corresponds to a fixed value of $n$, and the entries in that row are ordered with $k$ from $1$ to $n$.}
\label{tab:321}
\end{table}

\begin{table}[H]
\begin{tabular}{lllllllllllllllll}
  &   &    &    &    &    &     &     & 1   &     &      &     &      &      &      &      &      \\
  &   &    &    &    &    &     & 1   &     & 2   &      &     &      &      &      &      &      \\
  &   &    &    &    &    & 1   &     & 4   &     & 5    &     &      &      &      &      &      \\
  &   &    &    &    & 1  &     & 6   &     & 13  &      & 14  &      &      &      &      &      \\
  &   &    &    & 1  &    & 8   &     & 24  &     & 41   &     & 42   &      &      &      &      \\
  &   &    & 1  &    & 10 &     & 38  &     & 85  &      & 131 &      & 132  &      &      &      \\
  &   & 1  &    & 12 &    & 55  &     & 150 &     & 291  &     & 428  &      & 429  &      &      \\
  & 1 &    & 14 &    & 75 &     & 240 &     & 551 &      & 994 &      & 1429 &      & 1430 &      \\
1 &   & 16 &    & 98 &    & 359 &     & 945 &     & 1974 &     & 3424 &      & 4861 &      & 4862
\end{tabular}
\caption{(New sequence) The values $T(n,k)$ for $1\le n\le 9$ and $1\le k\le n$ such that $T(n,k)$ corresponds to the number of $312$-avoiding parking permutations of length $n$ with lucky set $I=[k]$ Here each row corresponds to a fixed value of $n$, and the entries in that row are ordered with $k$ from $1$ to $n$.}
\label{tab:312}
\end{table}

\begin{table}[H]
\begin{tabular}{lllllllllllllllll}
  &   &   &   &    &    &     &     & 1   &     &      &     &      &      &      &      &      \\
  &   &   &   &    &    &     & 1   &     & 2   &      &     &      &      &      &      &      \\
  &   &   &   &    &    & 1   &     & 3   &     & 5    &     &      &      &      &      &      \\
  &   &   &   &    & 1  &     & 4   &     & 9   &      & 14  &      &      &      &      &      \\
  &   &   &   & 1  &    & 5   &     & 14  &     & 28   &     & 42   &      &      &      &      \\
  &   &   & 1 &    & 6  &     & 20  &     & 48  &      & 90  &      & 132  &      &      &      \\
  &   & 1 &   & 7  &    & 27  &     & 75  &     & 165  &     & 297  &      & 429  &      &      \\
  & 1 &   & 8 &    & 35 &     & 110 &     & 275 &      & 572 &      & 1001 &      & 1430 &      \\
1 &   & 9 &   & 44 &    & 154 &     & 429 &     & 1001 &     & 2002 &      & 3432 &      & 4862
\end{tabular}
\caption{(Catalan's triangle (read by rows), OEIS \seqnum{A009766}) The values $T(n,k)$ for $1\leq n\leq 9$ and $1\leq k\leq n$ such that $T(n,k)$ corresponds to the number of $\sigma$-avoiding parking permutations of length $n$ with lucky set $I=[k]$ for $\sigma\in\{132, 213,  231\}$. Here each row corresponds to a fixed value of $n$, and the entries in that row are ordered with $k$ from $1$ to $n$. }
\label{tab:132-213-231}
\end{table}

\bibliographystyle{plainurl}
\bibliography{references.bib}

\end{document}